\documentclass[11pt,reqno]{amsart}
\usepackage[utf8x]{inputenc}

\usepackage{hyperref}
\hypersetup{
    colorlinks=true,       
    linkcolor=blue,          
    citecolor=blue,        
    filecolor=blue,      
    urlcolor=cyan  
}
\usepackage{latexsym,amssymb,amsmath,amsthm,bbm}
\usepackage{epsfig,mathrsfs}
\usepackage{mathtools}
\usepackage{mathabx}

\newcommand{\R}{{\ensuremath{\mathbb{R}}}}
\newcommand{\N}{{\ensuremath{\mathbb{N}}}}
\newcommand{\Z}{{\ensuremath{\mathbb{Z}}}}

\newcommand{\C}{{\ensuremath{\mathbb{C}}}}

\renewcommand{\P}{\ensuremath{\mathbb{P}}}
\renewcommand{\dj}{d\kern-0.4em\char"16\kern-0.1em}
\newcommand{\E}{\ensuremath{\mathbb{E}}}

\newcommand{\sig}{\ensuremath{\mathcal{F}}}

\newtheorem{Thm}{Theorem}[section]

\newtheorem{Lem}[Thm]{Lemma}
\newtheorem{Prop}[Thm]{Proposition}

\theoremstyle{definition}
\newtheorem{Rem}[Thm]{Remark}

\theoremstyle{definition}
\newtheorem{Ex}[Thm]{Example}

\theoremstyle{definition}

\theoremstyle{definition}

\begin{document}
\numberwithin{equation}{section}
\bibliographystyle{amsalpha}

\title[On subordinate random walks]{On subordinate random walks}
\begin{abstract}
In this article subordination of random walks in $\R^d$ is considered. We prove that subordination of random walks in the sense of \cite{BSC} yields the same process as subordination of L\' evy processes (in the sense of Bochner). Furthermore, we prove that appropriately scaled subordinate random walk converges to a multiple of a rotationally $2\alpha$-stable process if and only if the Laplace exponent of the corresponding subordinator varies regularly at zero with index $\alpha\in (0,1]$\,.
\end{abstract}

\author{Ante Mimica}
\address{University of Zagreb, Department of Mathematics, Bijeni\v cka cesta 30, 10000 Zagreb, Croatia}
\curraddr{}
\thanks{Research supported by Croatian Science Foundation under the project 3526}
\email{amimica@math.hr}



\maketitle

\allowdisplaybreaks[3]

\section{Introduction}

A subordinator $T=(T_t)_{t\geq 0}$ is a non-decreasing L\' evy process (i.e. a stochastic process having stationary and independent increments) defined on a probability space $(\Omega,\sig,\P)$. From the definition it follows that $T$ takes values in $[0,\infty)$ and the Laplace transform of $T_t$ is given by
\[
\E e^{-\lambda T_t}=e^{-t\phi(\lambda)},\quad\lambda>0\,, 
\]
where $\phi$ is called the Laplace exponent of $T$. It is of the following form (see \cite[III.1]{Be})
\begin{equation}\label{eq:LE}
\phi(\lambda)=b\lambda+\int_{(0,\infty)}(1-e^{-\lambda s})\mu(ds)\,,
\end{equation}
where $b\geq 0$ is the drift of the subordinator and $\mu$ is a measure on $(0,\infty)$ satisfying $\int_{(0,\infty)}(1\wedge s)\mu(ds)<\infty$ called the L\' evy measure of $T$\,.
It is known that $\phi$ is a Bernstein function, meaning that $\phi:(0,\infty)\rightarrow \R$ is a $C^\infty$-function satisfying
\[
    \phi(\lambda)\geq 0\quad \text{ for all }\quad \lambda>0
\]
and
\[
 (-1)^{n-1}\phi^{(n)}(\lambda)\geq 0\quad \text{ for all } \quad \lambda>0 \text{ and }\,\,\,\,n\in \N
\]
(here $\phi^{(n)}$ denotes the $n$-th derivative of $\phi$).
It is known that every Bernstein function $\phi$ satisfying $\lim\limits_{\lambda\downarrow 0}\phi(\lambda)=0$ has a unique representation given by \eqref{eq:LE} (see \cite[Theorem 3.2]{SSV}).

Let $Z=(Z_n)_{n\geq 0}$ be a random walk in $\R^d$ independent of $T$, that is $Z_0=0$ and $Z_n=\zeta_1+\ldots+\zeta_n$ for $n\geq 1$, where $(\zeta_n)_{n\geq 1}$ is a sequence of independent and identically distributed random vectors in $\R^d$ defined on the probability space $(\Omega,\sig,\P)$ idependent of $T$\,. 

We define subordinate random walk $X=(X_n)_{n\geq 0}$  as 
\begin{equation}\label{eq:dsub1}
X_0=0\quad \text{ and }\quad X_n=\sum_{k=1}^n \xi_k\,\,\text{ for }\,\,n\geq 1,
\end{equation}
  where $\xi_1,\xi_2,\ldots$ are independent and identically distributed  having the following distribution
\begin{equation}\label{eq:dsub2}
	\P(\xi_1\in B)=b\P(Z_1\in B)+\sum_{m=1}^\infty \int_{(0,\infty)}\frac{q^{m-1}t^m}{m!}e^{-qt}\mu(dt)\P(Z_m\in B),\quad B\subset \R^d \,\,\text{Borel}
\end{equation}
and  $q>0$ is such that $\phi(q)=q$.

\begin{Rem}
The choice of $q>0$ is made so that $\P(\xi_1\in \R^d)=1$:
\begin{align*}
\P(\xi_1\in \R^d)&=b+q^{-1}\int_{(0,\infty)} (e^{qt}-1)e^{-qt}\mu(dt)=\frac{\phi(q)}{q}=1.
\end{align*}
The condition  $\phi(q)=q$ is not a big restriction, since we can always consider normalization $\psi(\lambda)=q\phi(q)^{-1}\phi(\lambda)$ of $\phi$ that satisfies the latter condition.  
\end{Rem}

This procedure is known as discrete subordination introduced in \cite{BSC}. In this article we will embed subordinate random walks into a L\' evy process (actually compound Poisson processes) by using  time change by a Poisson process (see \cite[Section 1.4]{S} and Section \ref{sec:levy}). This gives us an opportunity to link discrete subordination with subordination of L\' evy processes (in the sense of Bochner, see \cite[Chapter 13]{SSV}). It turns out that, understood as L\' evy processes just described, discrete subordination and subordination of the random walk by the subordinator give rise to  the same process as the following proposition states. 
\begin{Prop}\label{prop:comm} Let $Z=(Z_n)_{n\geq 0}$ be a random walk in $\R^d$, $T=(T_t)_{t\geq 0}$ an independent nonconstant subordinator and let $X=(X_n)_{n\geq 0}$ be the corresponding discrete subordinate random walk defined by (\ref{eq:dsub1}) and (\ref{eq:dsub2}). If $N=(N_t)_{t\geq 0}$ is a Poisson process with intensity $q>0$ independent of $X$, $T$ and $Z$,  then the processes $\hat{X}=(\hat{X}_t)_{t\geq 0}$ and $\tilde{X}=(\tilde{X}_t)_{t\geq 0}$ defined by
\[
    \hat{X}_t=X_{N_t}\qquad \text{ and }\qquad \tilde{X}_t=Z_{N_{T_t}},\quad t\geq 0
\]
are compound Poisson processes having the same charcteristics. 
\end{Prop}

Let $Z=(Z_n)_{n\geq 0}$ be the simple symmetric  random walk in $\Z^d$, that is $Z_0=0$, $Z_n=\zeta_1+\ldots +\zeta_n$, for $n\geq 1$ and $(\zeta_n)_{n\geq 1}$ is a sequence of independent and identically distributed random vectors in 
$\Z^d$ such that $\P(\zeta_1=e_j)=\P(\zeta_1=-e_j)=\frac{1}{2d}$, where $e_j=(0,\ldots,0,1,0,\ldots,0)\in \Z^d$ has $1$ at the $j$-th coordinate, $1\leq j\leq d$\,. 

Consider a nonconstant subordinator $T=(T_t)_{t\geq 0}$ with the Laplace exponent $\phi$ and a sequence $(\xi_n)_{n\geq 1}$ of independent and identically distributed random variables with the distribution given by (\ref{eq:dsub2}) with $q=1$. Note that $\phi$ strictly increases (so that the inverse $\phi^{-1}$ exists).
We aim to consider convergence of the processes $X^{(n)}=(X^{(n)}_t)_{t\geq 0}$, $n\in \N$ defined by
\begin{equation}\label{eq:X-def}
 X_t^{(n)}=\sqrt{\phi^{-1}(n^{-1})}\sum_{k=1}^{N_{nt}}\xi_k,\quad n\geq 0,
\end{equation}
where $N=(N_t)_{t\geq 0}$ is a stochastic process that is either a Poisson proces with intensity $1$ independent of the sequence $(\xi_n)_{n\geq 1}$  or $N_t=\lfloor t\rfloor$, $t\geq 0$. In the former case (see Proposition \ref{prop:comm})
\begin{equation}\label{eq:process1}
    X_t^{(n)}\overset{d}{=}\sqrt{\phi^{-1}(n^{-1})}Z_{N_{T_{nt}}},
\end{equation}
while in the latter case 
\begin{equation}\label{eq:process2}
    X_t^{(n)}=\sqrt{\phi^{-1}(n^{-1})}\sum_{k=1}^{\lfloor nt\rfloor}\xi_k.
\end{equation}

Our aim is to investigate (weak) limits of the processes $X^{(n)}$, i.e. we will establish a type of a functional limit theorem by imposing conditions on the Laplace exponent $\phi$. Note that the processes defined by \eqref{eq:process1} are actually random sums and find their application in economics, financial mathematics and queuing theory(see \cite{EKM},\cite{G}).

The following tightness result in the space \[
D([0,\infty),\R^d):=\{f:[0,\infty)\rightarrow \R^d: f\text{ is right continuous with left limits}\}\] (with an appropriate topology) will be proved in Section \ref{sec:tight}\,. Condition that will ensure tightness for subordinate random walk is certain upper scaling condition of the inverse of the inverse of the Laplace exponent. 

\begin{Prop}\label{prop:tightness} Assume that the the Laplace exponent $\phi$ of a subordinator satisfies $\phi(1)=1$ and there exist $c>0$ and $\gamma>1$ such that  $\limsup\limits_{\lambda\downarrow 0}\frac{\phi^{-1}(\lambda x)}{\phi^{-1}(\lambda)}\leq cx^\gamma$ for all $x\geq 1$. 
Then the sequence $(X^{(n)})_{n\geq 1}$ defined by \eqref{eq:X-def} is tight in the space $D([0,\infty),\R^d)$\,.
\end{Prop}

Now we explore convergence of the sequence $(X^{(n)})_{n\geq 1}$. A measurable function $f\colon (0,\infty)\rightarrow (0,\infty)$ is said to vary regularly at $0$ with index $\rho\in \R$ if
\[
    \lim_{\lambda \downarrow 0}\frac{f(\lambda x)}{f(\lambda)}=x^\rho\quad \text{ for all }\quad  x>0\,.
\]
A L\' evy process (see Section \ref{sec:levy}) $W^{(\alpha)}=(W^{(\beta)}_t)_{t\geq 0}$ is called rotationally invariant $\beta$-stable process  in $\R^d$ ($\beta\in (0,2]$) if 
\[
    \E[e^{i\vartheta\cdot W_t^{(\beta)}}]=e^{-t|\vartheta|^\beta}\quad \text{ for all }\quad \vartheta\in \R^d\quad \text{ and }\quad t\geq 0\,.
\]

\begin{Thm}\label{thm:conv} Assume that the Laplace exponent $\phi$ of a subordinator satisfies $\phi(1)=1$.
The sequence $(X^{(n)})_{n\geq 1}$ converges in $D([0,\infty),\R^d)$ if and only if $\phi$ varies regularly at zero with index $\alpha\in (0,1]$. In this case the limit is $(2d)^{-\alpha}W^{(2\alpha)}$\,.
\end{Thm}

The following example is interesting, since it approximates Brownian motion by a compound Poisson process that has steps with infinite variance.

\begin{Ex}
Consider a Bernstein function 
\[
    \phi(\lambda)=c\lambda\log\left(1+\tfrac{1}{\lambda}\right),\quad c=\frac{1}{\log 2}\,.
\]
Since $\lim\limits_{\lambda \downarrow 0}\phi(\lambda)=0$, $\phi$ is the Laplace exponent of a subordinator. 
Note that $\phi$ varies regularly at $0$ with index $1$\,.
By Theorem \ref{thm:conv}, the corresponding subordinate random walk 
\[
    X_t^{(n)}=\sqrt{\phi^{-1}(n^{-1})}\sum_{i=1}^{\lfloor nt\rfloor} \xi_i
\]
converges to a multiple of Brownian motion by $(2d)^{-1}$. Let $d=1$. Unlike in Donsker's invariance principle, $\E \xi_1^2=+\infty$\,. Indeed,
\begin{align*}
    \E\xi_1^2&=\sum_{m=1}^\infty \E Z_m^2\int_{(0,\infty)} \frac{t^m}{m!}e^{-t}\mu(dt)=\sum_{m=1}^\infty m\int_{(0,\infty)} \frac{t^m}{m!}e^{-t}\mu(dt)\\
    &=\int_{(0,\infty)} \sum_{m=1}^\infty \frac{t^m}{(m-1)!}e^{-t}\mu(dt)=\int_{(0,\infty)} t\mu(dt)=\int_0^\infty \mu(t,\infty)\,dt.
\end{align*}
By \cite[pp. 230-231, table entry 27]{SSV}, the L\' evy measure of $\phi$ is given by
\[
    \mu(dt)=c\tfrac{1-e^{-t}(1+t)}{t^2}\,dt
\]
and so
\[
    \lim_{t\to\infty} \frac{\mu(t,\infty)}{\frac{1}{t}}=c.
\]
Therefore, there exists a constant $c_1>0$ so that \[
    \E \xi_1^2\geq c_1\int_1^\infty \frac{dt}{t}=+\infty\,.
\]
\end{Ex}


\section{L\' evy processes}\label{sec:levy}

A stochastic process $X=(X_t)_{t\geq 0}$ defined on a probability space $(\Omega,\sig,\P)$  taking values in $\R^d$ $(d\geq 1)$ is a L\' evy process if $X_0=0$, it has stationary and independent incremets and paths that are $\P$-a.s. right-continuous with left limits. It is well known that the characteristic function of $X_t$ is  
\[
\E e^{i\vartheta\cdot X_t}= e^{-t\psi(\vartheta)}\,,\quad t>0\,,\,\,\,\vartheta\in \R^d,
\]
where $\psi:\R^d\rightarrow \C$ is called the characteristic exponent of $X$ and has the following L\' evy-Khintchine representation
\begin{equation}
\psi(\vartheta)=i\beta\cdot \vartheta+\frac{1}{2}Q\vartheta\cdot \vartheta+\int_{\R^d\setminus\{0\}}(1-e^{iy\cdot \vartheta}+i y\cdot\vartheta 1_{\{|y|\leq 1\}})\nu(dy),\quad \vartheta\in \R^d.
\end{equation}
Here $\beta=(\beta_1,\ldots,\beta_d)\in \R^d$, $Q$ is a $d\times d$ positive semi-definite matrix and $\nu$ is the L\' evy measure, i.e. a measure on $\R^d\setminus\{0\}$ satisfying $\int_{\R^d\setminus\{0\}}(1\wedge |y|^2)\nu(dy)<\infty\,.$ If $\beta=0$ and $Q=0$, we will call $X$ a pure jump L\' evy process. A triplet $(\beta,Q,\nu)$ is called the L\' evy triplet of $X$ and every L\' evy process is uniquely determined by its L\' evy triplet (see \cite[Section 3.11]{S}). 

A particular example of a L\' evy process that will be used in this section is a compound Poisson process (see \cite[Section 1.4]{S}). It is a L\' evy process $X=(X_t)_{t\geq 0}$ with the characteristic exponent of the form
\[
\psi(\vartheta)=q\int_{\R^d}(1-e^{iy\cdot\vartheta})\eta(dy),\quad \vartheta\in \R^d,
\]
where $q>0$ and $\eta$ is a measure on $\R^d$ such that $\eta(\{0\})=0$\,. The corresponding L\' evy triplet is $(\beta,0,\nu)$, where 
\begin{equation}\label{eq:beta}
\beta=(\beta_1,\ldots,\beta_d),\,\,\,\beta_j=-q\int_{0<|y|\leq 1} y_j\eta(dy),\,\,j=1,\ldots,d
\end{equation}
and 
\begin{equation}\label{eq:levy}
\nu(B)=q\eta(B)\,,\quad B\subset \R^d\,\,\text{Borel}\,.
\end{equation}
In particular, for $\eta=\delta_{1}$, $X$ is the Poisson process with intensity $q$\,. 

The following construction of a compound Poisson process will be useful (see \cite[Theorem 1.4.2]{S}). Let $N=(N_t)_{t\geq 0}$ be the Poisson process with intensity $q>0$ and let $(\zeta_n)_{n\in \N}$ be a  sequence of independent random variables with law $\rho$ defined on a common probability space $(\Omega,\sig,\P)$ independent of $N$. Define a random walk $Z=(Z_n)_{n\geq 0}$ by $Z_0=0$ and $Z_n=\zeta_1+\ldots+\zeta_n$ for $n\geq 1$\,. Then the process $X=(X_t)_{t\geq 0}$ defined by $X_t=Z_{N_t}$ is a compound Poisson process. 

\proof[Proof of Proposition \ref{prop:comm}]
The process $\tilde{Z}=(\tilde{Z}_t)_{t\geq 0}$ defined by $\tilde{Z}_t=Z_{N_t}$ is a compound Poisson process with (see \eqref{eq:beta} and \eqref{eq:levy}) 
\[
    \beta_{\tilde{Z}}=-q\int_{0<|y|\leq 1}y\P(\zeta_1\in dy)\qquad \text{ and }\qquad \nu_{\tilde{Z}}(dy)=q\P(\zeta_1\in dy).
\]
By \cite[Theorem 30.1]{S}, $\tilde{X}$ is a L\' evy process with the L\' evy triplet $(\beta_{\tilde{X}}, 0, \nu_{\tilde{X}})$, where  
\begin{equation}\label{eq:comm1}
    \beta_{\tilde{X}}=b\beta_{\tilde{Z}}-\int_{(0,\infty)}\int_{0<|y|\leq 1}y\P(\tilde{Z}_s\in dy)\mu(ds)
\end{equation}
and
\begin{equation}\label{eq:comm2}
\nu_{\tilde{X}}(B)=qb\P(\zeta_1\in B)+\int_{(0,\infty)}\P(\tilde{Z_s}\in B)\mu(ds).
\end{equation}
On the other hand, $\hat{X}$ is a compound Poisson process and its L\' evy triplet  $(\beta_{\hat{X}}, 0, \nu_{\hat{X}})$ is, by (\ref{eq:beta}) and (\ref{eq:levy}), given by
\[
    \beta_{\hat{X}}=-q\int_{0<|y|\leq 1}y\P(\xi_1\in dy)\quad \text{and} \quad \nu_{\hat{X}}(dy)=q\P(\xi_1\in dy).
\]
Then (\ref{eq:dsub2}) yields
\begin{align*}
\beta_{\hat{X}}&=-qb\int_{0<|y|\leq 1}y\P(Z_1\in dy)-\int_{(0,\infty)}\int_{0<|y|\leq 1}\sum_{m=1}^\infty\frac{(qs)^m}{m!}e^{-qs}y\P(Z_m\in dy)\mu(ds)\\
&=-qb\int_{0<|y|\leq 1}y\P(\zeta_1\in dy)-\int_{(0,\infty)}\int_{0<|y|\leq 1} y\P(Z_{N_s}\in dy)\mu(ds)=\beta_{\tilde{X}}
\end{align*}
and, for a Borel set $B\subset \R^d\setminus\{0\}$, 
\begin{align*}
\nu_{\hat{X}}(B)&=qb\P(Z_1\in B)+\int_{(0,\infty)}\sum_{m=1}^\infty \frac{(qs)^m}{m!}e^{-qs}\P(Z_m\in B)\mu(ds)\\
&=qb\P(\zeta_1\in B)+\int_{(0,\infty)}\P(Z_{N_s}\in B)\mu(ds)=\nu_{\tilde{X}}(B)\,.
\end{align*}
Hence, $\hat{X}$ and $\tilde{X}$ are L\' evy processes with same L\' evy triplets.
\qed
\section{Tightness result}\label{sec:tight}
We recall that the process $X^{(n)}$ was defined by
\[
 X_t^{(n)}=\sqrt{\phi^{-1}(n^{-1})}\sum_{k=0}^{N_{nt}}\xi_k,\quad n\geq 0,
\]
where $T=(T_t)_{t\geq 0}$ is a subordinator with the Laplace exponent $\phi$, $(\xi_n)_{n\geq 1}$ is a sequence of independent and identically distributed random vectors with the distribution given by (\ref{eq:dsub2}) and $N=(N_t)_{t\geq 0}$ is either the Poisson process with intensity $1$ independent of the sequence $(\xi_n)_{n\geq 1}$ or $N_t=\lfloor t\rfloor$ for $t\geq 0$\,. 
The aim of this section is to prove that the sequence $(X^{(n)})_{n\geq 1}$ is tight in the Skorokhod space $D([0,\infty),\R^d)$ endowed with the Skorokhod topology. We refer the reader to  \cite[VI.1b]{JS} for a definition of the Skorokhod topology.

\begin{Lem}\label{lem:lmest}
Let $\phi:(0,\infty)\rightarrow (0,\infty)$ be a the Laplace exponent of a subordinator with the L\' evy measure $\mu$. 
\begin{itemize}
\item[(i)] Then \[
    \int_{(0,r]} t\mu(dt)\leq er\phi(r^{-1})\quad \text{ for all }\quad r>0
\]
and
\[
    \mu(t,\infty)\leq(1-e^{-1})^{-1}\phi(t^{-1})\quad \text{ for all }\quad t>0\,.
\]
\item[(ii)] For all $\lambda,x>0$ we have $\phi(\lambda x)\leq (x\vee 1)\phi(\lambda)$.
\item[(iii)] If $\phi$ is trictly increasing and  varies regularly at $0$ with index $\alpha>0$, then $\phi^{-1}$ varies regularly at $0$ with index $1/\alpha$.
\end{itemize}

\end{Lem}
\begin{Rem} Lemma \ref{lem:lmest} (iii) can be compared with \cite[Proposition 1.5.15]{BGT}, where asymptotic inverse and conjugacy are considered. Although this is a special case of the aforementioned result, we give simplified proof for strictly increasing functions $\phi$. 
\end{Rem}

\proof
    (i) Starting from \eqref{eq:LE} and using an elementary inequality \[1-e^{-x}\geq xe^{-x},\,x\geq 0\]it follows that 
    \[
        \phi(\lambda)\geq \int_{(0,{\lambda^{-1}}]}(1-e^{-\lambda t})\mu(dt)\geq \int_{(0,\lambda^{-1}]}\lambda te^{-\lambda t}\mu(dt)\geq e^{-1}\lambda\int_{(0,\lambda^{-1}]}t\mu(dt),
    \]
    yielding the first estimate. By \eqref{eq:LE},
    \[
        \phi(\lambda)\geq \int_{(\lambda^{-1},\infty)}(1-e^{-\lambda t})\mu(dt)\geq (1-e^{-1})\mu(\lambda^{-1},\infty)
    \]
    and the second estimate follows. 
    \\(ii)
    Let $\lambda>0$. If $x\geq 1$, then $1-e^{-tx}\leq x(1-e^{-t})$ for all $t\geq 0$, hence
    \begin{align*}
    \phi(\lambda x)&=b\lambda x+\int_{(0,\infty)}(1-e^{-\lambda x t})\mu(dt)\leq b\lambda x+\int_{(0,\infty)}x(1-e^{-\lambda t})\mu(dt)=x\phi(\lambda)\,.
    \end{align*}
    For $x\leq 1$  we use that $\phi$ is non-decreasing to conclude that $\phi(\lambda x)\leq \phi(\lambda)$. \\
    (iii) Let $x>0$. By regular variation of $\phi$, for any $0<\varepsilon <x^\alpha$ there exists $\delta>0$ such that 
    \[
        x^{\alpha}-\varepsilon\leq \frac{\phi(\lambda x)}{\phi(\lambda)}\leq x^{\alpha}+\varepsilon,\,\,\,0<\lambda<\delta.
    \]
    Since $\phi$ is strictly increasing and continuous, the last display implies
    \[
        \limsup_{\lambda\downarrow 0}\frac{\phi^{-1}(\lambda(x^\alpha-\varepsilon))}{\phi^{-1}(\lambda)}\leq x\leq \liminf_{\lambda\downarrow 0}\frac{\phi^{-1}(\lambda(x^\alpha+\varepsilon))}{\phi^{-1}(\lambda)}.
    \]
    Let $y>0$. By taking $x=(y+\varepsilon)^{1/\alpha}$ we get
    $
        \limsup\limits_{\lambda\downarrow 0}\frac{\phi^{-1}(\lambda y)}{\phi^{-1}(\lambda)}\leq (y+\varepsilon)^{1/\alpha}
    $ and, since $\varepsilon>0$ was arbitrary, this implies $ \limsup\limits_{\lambda\downarrow 0}\frac{\phi^{-1}(\lambda y)}{\phi^{-1}(\lambda)}\leq y^{1/\alpha}$. For $0<\varepsilon<y$ and $x=(y-\varepsilon)^{1/\alpha}$ we obtain 
    $\liminf\limits_{\lambda\downarrow 0}\frac{\phi^{-1}(\lambda y)}{\phi^{-1}(\lambda)}\geq (y-\varepsilon)^{1/\alpha}$ and so $\liminf\limits_{\lambda\downarrow 0}\frac{\phi^{-1}(\lambda y)}{\phi^{-1}(\lambda)}\geq y^{1/\alpha}$.
\qed

\begin{Lem} \label{lem:line-tightness}
Assume that $\phi$ is strictly increasing with $b=0$ and $\phi(1)=1$. There exists a constant $c_0>0$ such that for any $K,\beta,a>0$ and $n\in \N$ we have
\[
    \P(|\sum_{k\leq an}\xi_k|>\tfrac{K}{\sqrt{\phi^{-1}(n^{-1})}})\leq c_0 a(K^{-2-\beta}\tfrac{\phi^{-1}(n^{-1}))}{\phi^{-1}(K^{-\beta}n^{-1})}+K^{-\beta})\,.
\]
\end{Lem}
\proof
    Let $\{(Z_n^{(k)})_{n\geq 1}:k\in \N\}$ be a family of independent copies of the random walk $(Z_n)_{n\geq 1}$ defined by 
    $Z_0=0$, $Z_n=\zeta_1+\ldots +\zeta_n$, for $n\geq 1$,  where $(\zeta_n)_{n\geq 1}$ is a sequence of independent and identically distributed random vectors in $\Z^d$ such that $\P(\zeta_1=e_j)=\P(\zeta_1=-e_j)=\frac{1}{2d}$, $j=1,\ldots,d$. It follows from Chebyshev inequality that, for $r>0$ and $l,m_1,\ldots,m_l\in \N$ the following holds
    \[
        \P(|\sum_{k=1}^lZ_{m_k}^{(k)}|>r)\leq\tfrac{\E|\sum_{k=1}^lZ_{m_k}^{(k)}|^2}{r^2}= \tfrac{m_1+\ldots+m_l}{r^2}
    \]
    and this gives the following estimate
    \begin{equation}\label{eq:est-cheb2}
    \P(|\sum_{k=1}^lZ_{m_k}^{(k)}|>r)\leq \tfrac{m_1}{r^2}\wedge 1+\ldots+\tfrac{m_l}{r^2}\wedge 1.
    \end{equation}
By \eqref{eq:dsub2} with $q=1$ and \eqref{eq:est-cheb2}, 
\begin{align*}
\P(|\sum_{k\leq an}\xi_k|>&\tfrac{K}{\sqrt{\phi^{-1}(n^{-1})}})\\&=\sum_{m_1,\ldots,m_{\lfloor an\rfloor }=1}^\infty \prod_{j\leq an}\int_{(0,\infty)}\frac{t^{m_j}}{m_j!}e^{-t}\mu(dt)\P(|\sum_{k\leq an}Z_{m_k}^{(k)}|>\tfrac{K}{\sqrt{\phi^{-1}(n^{-1})}})\\
&\leq \sum_{i\leq an}\sum_{m_i=1}^\infty\int_{(0,\infty)}\frac{t^{m_i}}{m_i!}e^{-t}\mu(dt)(\tfrac{m_i\phi^{-1}(n^{-1})}{K^2}\wedge 1)\prod_{\substack{j\leq an\\j\not= i}}(\sum_{m_j=1}^\infty \int_{(0,\infty)}\tfrac{t^{m_j}}{m_j!}e^{-t}\mu(dt)) .
\end{align*}
Since $\sum\limits_{m=1}^\infty \int_{(0,\infty)}\frac{t^{m}}{m!}e^{-t}\mu(dt)=\int_{(0,\infty)}(e^t-1)e^{-t}\mu(dt)=\phi(1)=1$, the product term in the last display is equal to $1$ and so, for $\beta>0$, 
\begin{align*}
\P(|\sum_{k\leq an}&\xi_k|>\tfrac{K}{\sqrt{\phi^{-1}(n^{-1})}})\leq an\int_{(0,\infty)}\sum_{m=1}^\infty \tfrac{t^m}{m!}(\tfrac{m_i\phi^{-1}(n^{-1})}{K^2}\wedge 1)\mu(dt)\\
&\leq an(\int_{(0,\phi^{-1}(K^{-\beta}n^{-1})^{-1}]}t\mu(dt)K^{-2}\phi^{-1}(n^{-1})+\int_{(\phi^{-1}(K^{-\beta}n^{-1})^{-1},\infty)}\mu(dt)).
\end{align*}
Now we apply  Lemma \ref{lem:lmest} (i) 
to get 
\begin{align*}
\allowdisplaybreaks
\P(|\sum_{k\leq an}&\xi_k|>\tfrac{K}{\sqrt{\phi^{-1}(n^{-1})}})\\&\leq an(e\phi^{-1}(K^{-\beta}n^{-1})^{-1}K^{-\beta}n^{-1}K^{-2}\phi^{-1}(n^{-1})+\tfrac{e}{e-1}K^{-\beta}n^{-1})\\
&=aeK^{-2-\beta}\tfrac{\phi^{-1}(n^{-1})}{\phi^{-1}(K^{-\beta}n^{-1})}+\tfrac{ae}{e-1}K^{-\beta}.
\end{align*}
\qed

Recall that a random time $\tau$ is a stopping time with respect to the process $X^{(n)}$ if $\{\tau\leq t\}\in \sigma(X_s^{(n)}:s\leq t)$ for any $t>0$. In particular, if $N_t=\lfloor t\rfloor$, then $ \sigma(X_s^{(n)}:s\leq t)=\sigma(\xi_k:k\leq nt)$ and so
\begin{equation}\label{eq:stoppingtime}\{\lfloor n\tau\rfloor=m\}=\{\tfrac{m}{n}\leq \tau<\tfrac{m+1}{n}\}\in \sigma(\xi_k:k\leq m),\,\,\,m\in \N\cup\{0\}.\end{equation}
\begin{Lem}\label{lem:aux5}
Let $(\tau_n)_{n\geq 1}$ be a sequence of random times such that $\tau_n$ is a stopping time for the process $X^{(n)}$ for any $n\in \N$ and let $(h_n)$ be a sequence of non-negative numbers. There exists a sequence of functions $g_n:\R^d\rightarrow [0,\infty)$ satisfying $\lim\limits_{n\to\infty} g_n(\vartheta)=1$ for all $\vartheta\in \R^d$  such that 
\[
    E[e^{i\vartheta \cdot(X_{\tau_n+h_n}^{(n)}-X_{\tau_n}^{(n)})}]=\E[e^{i\vartheta\cdot X_{h_n}^{(n)}}]g_n(\vartheta),\quad \vartheta\in \R^d.
\]
\end{Lem}
\proof
If $N$ is a Poisson process, then $X^{(n)}$ is a L\' evy process and so the claim follows from the strong Markov property for L\' evy processes with $g_n\equiv 1$, $n\in \N$.

Let us consider now the case $N_t=\lfloor t\rfloor$. First we remark that 
\[
    Z_n:=N_{\tau_n+h_n}-N_{\tau_n}-N_{h_n}=\lfloor n(\tau_n+h_n)\rfloor -\lfloor n\tau_n\rfloor -\lfloor nh_n\rfloor\in \{0,1\}\,\,\,\P-\text{a.s.},\,\,n\in \N.
\]
For  $\vartheta\in \R^d$ we calculate
\begin{align*}
\E[&e^{i\vartheta \cdot (X_{\tau_n+h_n}^{(n)}-X_{\tau_n}^{(n)})}]\\&=\E[\exp\{i\sqrt{\phi^{-1}(n^{-1})}\sum_{k=\lfloor n\tau_n\rfloor +1}^{\lfloor n\tau_n\rfloor +\lfloor nh_n\rfloor}\vartheta\cdot \xi_k\}e^{i\sqrt{\phi^{-1}(n^{-1}))}Z_n\vartheta\cdot \xi_{\lfloor n\tau_n\rfloor +\lfloor nh_n\rfloor +1}}]\\
&=\sum_{m=0}^\infty \E[\exp\{i\sqrt{\phi^{-1}(n^{-1})}\sum_{k=m +1}^{m +\lfloor nh_n\rfloor}\vartheta\cdot \xi_k\}e^{i\sqrt{\phi^{-1}(n^{-1}))}Z_n\vartheta\cdot \xi_{m +\lfloor nh_n\rfloor +1}};\lfloor n\tau_n\rfloor=m].
\end{align*}
Now we use that $(\xi_k)_{k\geq 1}$ are independent and identically distributed and \eqref{eq:stoppingtime} to deduce
\begin{align*}
\E[&e^{i\vartheta \cdot (X_{\tau_n+h_n}^{(n)}-X_{\tau_n}^{(n)})}]\\
&=\sum_{m=0}^\infty \E[\exp\{i\sqrt{\phi^{-1}(n^{-1})}\sum_{k=m +1}^{m +\lfloor nh_n\rfloor}\vartheta\cdot \xi_k\}]\E[e^{i\sqrt{\phi^{-1}(n^{-1}))}Z_n\vartheta\cdot \xi_{m +\lfloor nh_n\rfloor +1}}]\P(\lfloor n\tau_n\rfloor=m)\\
&=\E[\exp\{i\sqrt{\phi^{-1}(n^{-1})}\sum_{k=1}^{\lfloor nh_n\rfloor}\vartheta\cdot \xi_k\}]\E[e^{i\sqrt{\phi^{-1}(n^{-1}))}Z_n\vartheta\cdot \xi_1}].
\end{align*}
By dominated convergence theorem,
\[
    g_n(\vartheta):=\E[e^{i\sqrt{\phi^{-1}(n^{-1}))}Z_n\vartheta\cdot \xi_1}],\,\,\vartheta\in \R^d,\,n\in \N
\]
satisfies $\lim\limits_{n\to\infty} g_n(\vartheta)=1$ for all $\vartheta\in \R^d$.

\qed

\begin{Lem}\label{lem:aux6}
For any $\vartheta=(\vartheta_1,\ldots,\vartheta_d)\in \R^d$ and $t>0$ we have 
\[
    \E[e^{i\vartheta\cdot X_t^{(n)}}]=\E[(\phi(1)-\phi(1-
    \tfrac{\cos{\sqrt{\phi^{-1}(n^{-1})}\vartheta_1}+\ldots+\cos{\sqrt{\phi^{-1}(n^{-1})}\vartheta_d}}{d}
    )^{N_t}].
\]
\end{Lem}
\proof
Since $(\xi_k)_{k\geq 1}$ is independent,  identically distributed and independent of $N$, we get
\[
    \E[e^{i\vartheta\cdot X_t^{(n)}}]=\E[\E[e^{i\sqrt{\phi^{-1}(n^{-1})\vartheta\cdot \xi_1}}]^{N_t}].
\]
It is enough to note that 
\begin{align*}
\E[e^{i\sqrt{\phi^{-1}(n^{-1})}\vartheta\cdot\xi_1}]&=\sum_{m=1}^\infty \E[e^{i\sqrt{\phi^{-1}(n^{-1})}\vartheta\cdot Z_m}]\int_{(0,\infty)}\tfrac{t^m}{m!}e^{-t}\mu(dt)\\
&=\sum_{m=1}^\infty \int_{(0,\infty)}\tfrac{(t\E[e^{i\sqrt{\phi^{-1}(n^{-1})}\vartheta\cdot \zeta_1}])^m}{m!}\mu(dt)\\
&=\int_{(0,\infty)}(e^{t\E[e^{i\sqrt{\phi^{-1}(n^{-1})}\vartheta\cdot \zeta_1}]}-1)e^{-t}\mu(dt)\\
&=\phi(1)-\phi(1-\E[e^{i\sqrt{\phi^{-1}(n^{-1})}\vartheta\cdot \zeta_1}]).
\end{align*}
Now it is enough to note that 
\begin{equation*}\label{eq:chf}
    \E[e^{i\sqrt{\phi^{-1}(n^{-1})}\vartheta\cdot \zeta_1}]=\tfrac{\cos{\sqrt{\phi^{-1}(n^{-1})}\vartheta_1}+\ldots+\cos{\sqrt{\phi^{-1}(n^{-1})}\vartheta_d}}{d}.
\end{equation*}
\qed

\proof[Proof of Proposition \ref{prop:tightness}]

To prove tightness we use Aldous' criterion (see \cite[Theorem 1]{Al}), that is we show that for a sequence of bounded stopping times $(\tau_n)_{n\geq 1}$  (with respect to the natural filtration of $X^{(n)}$ for any $n\in \N$) and a sequence  $(h_n)_{n\geq 1}$ of positive numbers converging to $0$ the following holds
\begin{equation}\label{eq:prob-conv}
    Y_n:=X_{\tau_n+h_n}^{(n)}-X_{\tau_n}^{(n)} \text{ converges to } 0 \text{ in probability as } n\to\infty
\end{equation}
and, for any $t>0$ the following tightness of $(X_t^{(n)})_{n\geq 1}$: for any $\varepsilon>0$ there exists $K>0$ such that 
\begin{equation}\label{eq:proof-tightness}
\limsup_{n\to\infty}\P(|X_t^{(n)}|>K)<\varepsilon.
\end{equation}
In order to prove \eqref{eq:prob-conv}, it is enough to prove convergence in distribution to $0$. By Lemma \ref{lem:aux5} and Lemma \ref{lem:aux6}, for any $\vartheta\in \R^d$ we have
\begin{align}\nonumber
    \lim_{n\to\infty}\E[e^{i\vartheta\cdot Y_n}]&=\lim_{n\to\infty}\E[e^{i\vartheta\cdot X_{h_n}^{(n)}}]\\&=\lim_{n\to\infty}\E[(1-\phi(1-\tfrac{\cos{\sqrt{\phi^{-1}(n^{-1})}\vartheta_1}+\ldots+\cos{\sqrt{\phi^{-1}(n^{-1})}\vartheta_d}}{d}))^{N_{nh_n}}].\label{eq:ald11}
\end{align}
Assume first that $N$ is the Poisson process with intensity $1$. By using \eqref{eq:ald11} and formula
\begin{equation}\label{eq:formulapp}
\E[\theta^{N_t}]=e^{t(\theta-1)}, \quad \theta,t> 0
\end{equation}
we obtain
\begin{align*}
\lim_{n\to\infty}\E[e^{i\vartheta\cdot Y_n}]=\lim_{n\to\infty}e^{-nh_n\phi(1-\frac{\cos{\sqrt{\phi^{-1}(n^{-1})}\vartheta_1}+\ldots+\cos{\sqrt{\phi^{-1}(n^{-1})}\vartheta_d}}{d})},\quad \vartheta\in \R^d.
\end{align*}
To evaluate the limit we use Lemma \ref{lem:lmest} (i)
 to get
\begin{align}
0&\leq nh_n\phi (1-\tfrac{\cos{\sqrt{\phi^{-1}(n^{-1})}\vartheta_1}+\ldots+\cos{\sqrt{\phi^{-1}(n^{-1})}\vartheta_d}}{d})\nonumber \\&\leq nh_n\phi(\tfrac{\phi^{-1}(n^{-1})|\vartheta|^2}{2d})
\leq nh_n (\tfrac{|\vartheta|^2}{2d}\vee 1)\phi(\phi^{-1}(n^{-1}))=h_n(\tfrac{|\vartheta|^2}{2d}\vee 1).\label{eq:phi-tmp121}
\end{align}
Then  $\lim\limits_{n\to\infty}\E[e^{i\vartheta\cdot Y_n}]=1$ 
and so convergence in distribution to $0$ by continuity theorem.
Further, by Lemma \ref{lem:line-tightness} it follows that, for any $\beta,t>0$ and $K\geq 1$, 
\begin{align*}
\P(|X_t^{(n)}|>K)&=\P(|\sum_{k\leq N_{nt}}\xi_k|>\tfrac{K}{\sqrt{\phi^{-1}(n^{-1}))}})\\
&=\sum_{m=1}^\infty \P(|\sum_{k\leq \tfrac{m}{n}n}\xi_k|>\tfrac{K}{\sqrt{\phi^{-1}(n^{-1})}})\P(N_{nt}=m)\\
&\leq \sum_{m=1}^\infty c_0\tfrac{m}{n}(K^{-2-\beta}\tfrac{\phi^{-1}(n^{-1})}{\phi^{-1}(K^{-\beta}n^{-1})}+K^{-\beta})\tfrac{(nt)^m}{m!}e^{-nt}\\
&=c_0(K^{-2-\beta}\tfrac{\phi^{-1}(n^{-1})}{\phi^{-1}(K^{-\beta}n^{-1})}+K^{-\beta})t\sum_{m=1}^\infty \tfrac{(nt)^{m-1}}{(m-1)!}e^{-nt}\\
&=c_0t(K^{-2-\beta}\tfrac{\phi^{-1}(n^{-1})}{\phi^{-1}(K^{-\beta}n^{-1})}+K^{-\beta}).
\end{align*}
Therefore
\[
    \limsup_{n\to\infty}\P(|X_t^{(n)}|>K)\leq c_1t(K^{-2-\beta+\beta\gamma}+K^{-\beta})
\]    
for some constant $c_1>0$
 and so, by choosing $\beta>0$ small enough so that 
 \begin{equation}\label{eq:beta-choice}
-2-\beta+\beta\gamma<0,\quad \text{i.e.} \quad \beta<\frac{2}{\gamma-1}
\end{equation}
 and then $K\geq 1$ large enough, we get \eqref{eq:proof-tightness}, finishing the proof of tightness in this case.
 
Let $N_t=\lfloor t\rfloor$. From \eqref{eq:ald11}  we obtain
\begin{align*}
\lim_{n\to\infty}\E[e^{i\vartheta\cdot Y_n}]&=\lim_{n\to\infty} e^{-\lfloor nh_n\rfloor \ln(1-\phi(1-\tfrac{\cos{\sqrt{\phi^{-1}(n^{-1})}\vartheta_1}+\ldots+\cos{\sqrt{\phi^{-1}(n^{-1})}\vartheta_d}}{d}))}\\
&=\lim_{n\to\infty}e^{\lfloor nh_n\rfloor\phi(1-\frac{\cos{\sqrt{\phi^{-1}(n^{-1})}\vartheta_1}+\ldots+\cos{\sqrt{\phi^{-1}(n^{-1})}}}{d})}=1
\end{align*}
by \eqref{eq:phi-tmp121}, showing the required convergence in distribution.
Further,
by Lemma \ref{lem:line-tightness}, for any $t,\beta >0$ and $K\geq 1$
\begin{align*}
\limsup_{n\to\infty}\P(|X_t^{(n)}|>K)&=\limsup_{n\to\infty}\P(|\sum_{k\leq nt}\xi_k|>\tfrac{K}{\sqrt{\phi^{-1}(n^{-1})}})\\
&\leq c_0t\limsup_{n\to\infty}(K^{-2-\beta}\tfrac{\phi^{-1}(n^{-1})}{\phi^{-1}(K^{-\beta}n^{-1})}+K^{-\beta})\leq c_1t(K^{-2-\beta+\beta\gamma}+K^{-\beta}). 
\end{align*}
Choosing $\beta>0$ as in \eqref{eq:beta-choice}
  we may choose $K\geq 1$ large enough so that 
\eqref{eq:proof-tightness} holds.

\qed

\section{Weak convergence}

We start with a few auxilliary results. 
\begin{Lem}\label{lem:ind-reg-var}
Let $(a_n)_{n\geq 1}$ and $(b_n)_{n\geq 1}$ be sequences of non-negative numbers such that $\lim\limits_{n\to\infty } a_n =a\in (0,\infty)$, $\lim\limits_{n\to\infty} b_n=0$ and let $f:(0,\infty)\rightarrow (0,\infty)$ be a monotone function that varies regularly at $0$ with index $\rho\in \R$\,. Then
\[
    \lim_{n\to\infty}\frac{f(a_nb_n)}{f(b_n)}=a^\rho.
\]    
\end{Lem}
\proof
    Assume that $f$ is non-decreasing. For any $0<\varepsilon<a$ and $n\in \N$ large enough we have
    \begin{align*}
    \frac{f((a-\varepsilon)b_n)}{f(b_n)}\leq 
    \frac{f(a_nb_n)}{f(b_n)}\leq 
    \frac{f((a+\varepsilon)b_n)}{f(b_n)}.
    \end{align*}
This implies 
\[
    (a-\varepsilon)^\rho\leq \liminf_{n\to\infty} \frac{f(a_nb_n)}{f(b_n)}\leq \limsup_{n\to\infty} \frac{f(a_nb_n)}{f(b_n)}\leq (a+\varepsilon)^\rho.
\]
Letting $\varepsilon\to 0$ we obtain
\[
    \lim_{n\to\infty}\frac{f(a_nb_n)}{f(b_n)}=a^\rho.
\]
If $f$ is non-increasing, the proof is similar. 
\qed

\begin{Lem}\label{lem:aux1}
Let $\phi:(0,\infty)\rightarrow (0,\infty)$ be a Bernstein function.
\begin{itemize}
\item[(i)] If $\phi$ varies regularly at $0$ with index $\rho\in\R$, then $\rho\in [0,1]$\,.
\item[(ii)] Let $(a_n)_{n\geq 1}$ be a sequence of positive numbers such that $a_n\leq 1$ for all $n\in \N$ and $\lim\limits_{n\to\infty}a_n=1$. 
Then 
\[
    \lim_{n\to\infty}\frac{\phi(a_nx)}{\phi(x)}=1\text{ uniformly in } x>0.
\] 
\end{itemize}
 
\end{Lem}
\proof
    (i) Let $x\geq 1$\,. By Lemma \ref{lem:lmest} (ii),
 it follows that 
    \[
        x\geq \lim_{\lambda\downarrow 0}\frac{\phi(\lambda x)}{\phi(\lambda)}=x^\rho,
    \]
    hence $\rho\leq 1$. On the other hand, since $\phi$ is non-decreasing,
    \[
        1\leq \lim_{\lambda\downarrow 0}\frac{\phi(\lambda x)}{\phi(\lambda)}=x^\rho
    \]
and so $\rho\geq 0$\,.\\
(ii)  
By Lemma \ref{lem:lmest} (ii), for any $x>0$, $\phi(x)=\phi(a_nx a_n^{-1})\leq a_n^{-1}\phi(a_n x)$ and so
\[
    a_n\leq \frac{\phi(a_nx)}{\phi(x)}\leq 1
\]
yielding the claim. 
\qed

The following sufficient result for regular variation is taken from \cite{Fel}\,.

\begin{Lem}\label{lem:feller}
Let $(\lambda_n)_{n\geq 1}$ and $(a_n)_{n\geq 1}$ be sequences of positive numbers such that 
\[
    \lim_{n\to\infty} \frac{a_n}{a_{n+1}}=1\quad \text{ and }\quad \lim_{n\to \infty}\lambda_n=0\,.
\]
If $f:(0,\infty)\rightarrow (0,\infty)$ is monotone, 
\[
    g(y)=\lim_{n\to\infty}a_n f(\lambda_ny)\in [0,\infty]
\]
exists on a dense subset of $(0,\infty)$ and it is  finite and positive on some interval, then $f$ varies regularly at $0$ with index $\alpha\in \R$\,.
\end{Lem}
\proof
    See \cite[Lemma VIII.8.3]{Fel}\,.
\qed

\proof[Proof of Theorem \ref{thm:conv}]
 Assume that $X^{(n)}$ converges to $X=(X_t)_{t\geq 0}$ in $D([0,\infty),\R^d)$. Then $X_t^{(n)}$ converges in distribution to $X_t$ for any $t>0$ and, by 
 Lemma \ref{lem:aux6} the following limit exists
 \begin{equation}\label{eq:tmp103}
     \lim_{n\to\infty}\E[e^{iX_t^{(n)}\cdot \vartheta}]=\lim_{n\to\infty}\E[(1-\phi(1-\tfrac{\cos{\sqrt{\phi^{-1}(n^{-1})}\vartheta_1}+\ldots+\cos{\sqrt{\phi^{-1}(n^{-1})}\vartheta_d}}{d}))^{N_{nt}}]
 \end{equation}
 for any $\vartheta=(\vartheta_1,\ldots,\vartheta_d)\in \R^d$\,.
If $(N_t)_{t\geq 0}$ is the Poisson process with intensity $1$
by \eqref{eq:formulapp} and \eqref{eq:tmp103} 
the following limit exists
\begin{equation}\label{eq:tmp207}
    \lim_{n\to\infty}\E[e^{iX_t^{(n)}\cdot \vartheta}]=\lim_{n\to\infty}e^{-nt\phi\left(1-\frac{\cos{\sqrt{\phi^{-1}(n^{-1})}\vartheta_1}+\ldots+\cos{\sqrt{\phi^{-1}(n^{-1})}\vartheta_d}}{d}\right)}.
\end{equation}  
For $N_t=\lfloor t\rfloor$, we see that the following limit exists
\begin{align}\nonumber
\lim_{n\to\infty}\E[e^{iX_t^{(n)}\cdot \vartheta}]&=\lim_{n\to\infty} e^{-\lfloor nt\rfloor\log\left(1-\phi\left(1-\frac{\cos{\sqrt{\phi^{-1}(n^{-1})}\vartheta_1}+\ldots+\cos{\sqrt{\phi^{-1}(n^{-1})}\vartheta_d}}{d}\right)\right)}\\&=\lim_{n\to\infty}e^{-tn\phi\left(1-\frac{\cos{\sqrt{\phi^{-1}(n^{-1})}\vartheta_1}+\ldots+\cos{\sqrt{\phi^{-1}(n^{-1})}\vartheta_d}}{d}\right)}.\label{eq:tmp208}
\end{align}
Therefore, in both cases 
\[
    \lim_{n\to\infty}n\phi(1-\tfrac{\cos{\sqrt{\phi^{-1}(n^{-1})}\vartheta_1}+\ldots+\cos{\sqrt{\phi^{-1}(n^{-1})}\vartheta_d}}{d})
\]
exists.

Define  
\[
a_n=\frac{1-\frac{\cos{\sqrt{\phi^{-1}(n^{-1})}\vartheta_1}+\ldots+\cos{\sqrt{\phi^{-1}(n^{-1})}\vartheta_d}}{d}}{\frac{\phi^{-1}(n^{-1})|\vartheta|^2}{2d}}\quad \text{ and }\quad b_n=\phi^{-1}(n^{-1})\tfrac{|\vartheta|^2}{2d}.
\]
Since $a_n\leq 1$ for all $n\in \N$, $\lim\limits_{n\to\infty}a_n=1$ and $\lim\limits_{n\to\infty} b_n=0$, Lemma \ref{lem:aux1} implies that 
\begin{equation}
\lim_{n\to\infty}\frac{\phi(1-\tfrac{\cos{\sqrt{\phi^{-1}(n^{-1})}\vartheta_1}+\ldots+\cos{\sqrt{\phi^{-1}(n^{-1})}\vartheta_d}}{d})}{\phi(\phi^{-1}(n^{-1})\frac{|\vartheta|^2}{2d})}=\lim_{n\to\infty}\frac{\phi(a_n b_n)}{\phi(b_n)}=1;
\end{equation}
hence, in both cases, the limit
\begin{equation}\label{eq:tmp102}
    \lim_{n\to\infty} n\phi(\phi^{-1}(n^{-1})\tfrac{|\vartheta|^2}{2d})=\lim_{n\to\infty}n\phi\left(1-\tfrac{\cos{\sqrt{\phi^{-1}(n^{-1})}\vartheta_1}+\ldots+\cos{\sqrt{\phi^{-1}(n^{-1})}\vartheta_d}}{d}\right)
\end{equation}
exists. 
By Lemma \ref{lem:feller} (with $\lambda_n=\phi^{-1}(n^{-1})$ and $a_n=n$) we conclude that $\phi$ varies regularly at $0$ with index $\alpha\in \R$. Lemma \ref{lem:ind-reg-var} ensures that $0\leq \alpha\leq 1$. Note that \eqref{eq:tmp103},  \eqref{eq:tmp102} and \eqref{eq:tmp207} (or \eqref{eq:tmp208}) yield
\begin{equation}\label{eq:char_fn}
    \E[e^{i\vartheta\cdot X_t}]=\lim_{n\to\infty}\E[e^{i\vartheta\cdot X_t^{(n)}}]
    =\lim_{n\to\infty} e^{-tn\phi(\phi^{-1}(n^{-1})\frac{|\vartheta|^2}{2d})}
    =\lim_{n\to\infty} e^{-t\frac{\phi(\phi^{-1}(n^{-1})\frac{|\vartheta|^2}{2d})}{\phi(\phi^{-1}(n^{-1}))}}\,.
\end{equation}
Hence, if $\alpha\in (0,1]$
\[
    \E[e^{i\vartheta\cdot X_t}]=e^{-t(2d)^{-\alpha}|\vartheta|^{2\alpha}}\,,
\]
and so $X=(2d)^{-\alpha}W^{(2\alpha)}$ is a multiple of a rotationally invariant $2\alpha$-stable process. 
On the other hand, if $\alpha=0$ we get 
\[
    \E[e^{i\vartheta\cdot X_t}]=e^{-t}\,,
\]
which is impossible, since for $\vartheta\to 0$ we get a contradiction: $1=e^{-t}$ for all $t\geq 0$\,. 

Assume that $\phi$ varies regularly at $0$ with index $\alpha\in (0,1]$. By Lemma \ref{lem:lmest} (iii), $\phi^{-1}$ varies at $0$ with index $1/\alpha$ and so 
\[
    \limsup_{\lambda\downarrow 0}\frac{\phi^{-1}(\lambda x)}{\phi^{-1}(\lambda)}\leq x^{1/\alpha},\quad x\geq 1.
\]
Note that we can always replace $1/\alpha$ in the last display by any $\gamma>1/\alpha\geq 1$. Hence we may use Proposition \ref{prop:tightness} to obtain tightness 
 of the sequence $(X^{(n)})_{n\geq 1}$ in $D[(0,\infty),\R^d)$.
 To prove convergence, it is enough to prove convergence of finite-dimensional distributions (see \cite[Theorem 16.10]{Kal}). 
By Lemma \ref{lem:aux5}, for any $0\leq s<t$ (with $\tau_n=s$ i $h_n=t-s$) we get 
\begin{align*}
    \lim_{n\to\infty}\E[e^{i\vartheta\cdot (X_t^{(n)}-X_s^{(n)})}]=\lim_{n\to\infty}\E[e^{i\vartheta\cdot X_{t-s}^{(n)}}]\,.
\end{align*}
Hence, by \eqref{eq:char_fn}, 
\begin{equation}\label{eq:incr-chf}
    \lim_{n\to\infty}\E[e^{i\vartheta\cdot (X_t^{(n)}-X_s^{(n)})}]=e^{-(t-s)(2d)^{-\alpha} |\vartheta|^{2\alpha}}\,.
\end{equation}

	Let $j\in \N$, $0=t_0\leq t_1<t_2<\ldots<t_j$  and $\vartheta^{(1)},\vartheta^{(2)},\ldots,\vartheta^{(j)}\in \R^d$. Using independent increments property of L\' evy processes (if $N$ is Poisson process) or independence (if $N_t=\lfloor t\rfloor$) and (\ref{eq:incr-chf}) we get
	\begin{align*}
		\lim_{n\to\infty}\E\exp{\{i\sum_{k=1}^j \vartheta^{(k)}\cdot X_{t_k}^{(n)}\}}&=\lim_{n\to\infty}\E\exp{\{i\sum_{k=1}^j \sum_{l=1}^k\vartheta^{(k)}\cdot (X_{t_l}^{(n)}-X_{t_{l-1}}^{(n)})\}}\\
		&=\lim_{n\to\infty}\E\exp{\{i\sum_{l=1}^j \sum_{k=l}^j\vartheta^{(k)}\cdot (X_{t_l}^{(n)}-X_{t_{l-1}}^{(n)})\}}\\
		&=\prod_{l=1}^j\lim_{n\to\infty}\E\exp{\{i \sum_{k=l}^j\vartheta^{(k)}\cdot (X_{t_l}^{(n)}-X_{t_{l-1}}^{(n)})\}}\\
		&=\prod_{l=1}^j \exp{\{-(t_l-t_{l-1})(2d)^{-\alpha}|\sum_{k=l}^j\vartheta^{(k)}|^{2\alpha}\}}\\
		&=\E\exp{\{i\sum_{k=1}^j \vartheta^{(k)}\cdot X_{t_k}\}},
	\end{align*}
	where the last line is obtained by using independent increments property of a L\' evy process $X=(2d)^{-\alpha}W^{(\alpha)}$. Now it is enough to apply continuity theorem for characteristic functions to obtain  convergence in distrubution of $(X_{t_1}^{(n)},\ldots,X_{t_j}^{(n)})$ to $(X_{t_1},\ldots,X_{t_j})$. 
\qed

\bibliography{myrefs}
\end{document}